\documentclass[11pt]{amsart}
\usepackage{amssymb,amsmath,amsthm}
\usepackage{mathrsfs,dsfont,a4wide}

\usepackage{latexsym}
\usepackage{amsfonts}
\usepackage{amsmath}
\usepackage{amssymb}
\usepackage{graphicx}
\numberwithin{equation}{section}
\newtheorem{theorem}{Theorem}[section]

\newtheorem{definition}[theorem]{Definition}

\newtheorem{lemma}[theorem]{Lemma}

\newtheorem{remark}[theorem]{Remark}

\newtheorem*{acknowledgement}{Acknowledgement}

\begin{document}
\def\Re{\mathop{\rm Re}\,}
\def\Im{\mathop{\rm Im}\,}
\def\dom{\mathop{\rm dom}\,}
\def\dist{\mathop{\rm dist}}
\def\grad{\mathop{\rm grad}}


\author[G. Alessandrini]
{Giovanni Alessandrini}
\address[G. Alessandrini]{Dipartimento di Matematica e Geoscienze, Universit\`a  di Trieste, Via Valerio 12/b, 34100 Trieste, Italia} 
\email[G. Alessandrini]{alessang@units.it}
\author[V. Nesi]
{Vincenzo Nesi}
\address[V.  Nesi]{Dipartimento di Matematica ``G. Castelnuovo", Sapienza, Universit\'a di Roma,
Piazzale A. Moro 2, 00185 Roma, Italy}
\email[V. Nesi]{nesi@mat.uniroma1.it}

\title[Quantitative estimates on Jacobians]{Quantitative estimates on Jacobians\\for hybrid inverse problems}

\begin{abstract} 
 We consider  $\sigma$-harmonic mappings, that is mappings $U$ whose components $u_i$ solve a divergence structure elliptic equation ${\rm div} (\sigma \nabla u_i)=0$, for $i=1,\ldots,n $. We investigate whether, with suitably prescribed Dirichlet data, the Jacobian determinant can be bounded away from zero. Results of this sort are required in the treatment of the so-called hybrid inverse problems, and also in the field of homogenization studying bounds for the effective properties of composite materials.

\textit{2000 AMS Mathematics Classification Numbers: 30C62, 35J55}

\end{abstract}

\maketitle{}

\begin{flushright}
\textit{In memoria di Alfredo.
}\end{flushright}

\section{Introduction}
\hspace{0.7 cm}

\noindent
The appearance of coupled physics methods has provoked a sharp change of perspective in inverse boundary problems. The simultaneous use of different physical modalities to interrogate, through exterior measurements, a body whose interior parameters are unknown has enabled to single out interior functionals which carry
useful, and possibly stable, information on the parameters of interest. Such methods are also known under the name of ``hybrid inverse problems''. Notable examples are the coupling of Magnetic Resonance with Electrical Impedance Tomography \cite{seo-woo}, Ultrasound and Electrical Impedance Tomography \cite{ammari-capde}, Magnetic Resonance and Elastography \cite{honda-mclaugh}. To fix ideas, let us focus on Ultrasound Modulated Electrical Impedance Tomography. In EIT the goal is to determine the, possibly anisotropic, electrical conductivity $\sigma=\{\sigma_{ij}\}$ of a body $\Omega$ by repeated boundary measurements of voltage $u|_{\partial \Omega}$ and current distribution $\sigma \nabla u \cdot \nu$ with $u$ solving the elliptic PDE
\begin{equation}\label{basicpde}
{\rm div} (\sigma \nabla u)=0,~~\hbox{in $\Omega$}.
\end{equation}
As is well known \cite{mandache}, the stability is very weak and, in fact, in the anisotropic case, also non-uniqueness occurs \cite{kohn-vogelius}. By combining electrical measurements with ultrasound measurements it is possibile to focus on a tiny spot near any point $x\in \Omega$ and it has been shown by Ammari et al.\cite{ammari-capde} that one can detect the localized energy
\begin{equation}
H(x)=\sigma \nabla u\cdot \nabla u (x).
\end{equation}
If one repeats the experiments with different boundary voltages, it is possible to extract the functionals
\begin{equation}
H_{ij}(x)=\sigma \nabla u_i\cdot \nabla u_j (x).
\end{equation}
where $u_1,\hdots u_n$ is an array of different solutions to \eqref{basicpde}.
In Monard and Bal \cite{bal-monard1, bal-monard2}, it is shown how, from such functionals, one may obtain the conductivity $\sigma$ in a satisfactory stable fashion. The crucial point, however, is to be able to set up an array of boundary data $\phi_1\hdots,\phi_n$  and corresponding solutions $u_1,\hdots u_n$ in such a way that the functionals $H_{ij}$ are non degenerate.

\noindent
In other words, calling $U:\Omega \to \mathbb R^n$,  the mapping $U=(u_1,\hdots, u_n)$, which we shall designate ``$\sigma$-harmonic'' mapping, it is required that the Jacobian determinant
\begin{equation*}
\det DU
\end{equation*}
does not vanish. And, furthermore, for the purpose of stability, a quantitative lower bound would be needed.

\noindent
This is the main question that we wish to address in this note,  which essentially stays behind all coupled physics problems mentioned above, and other inverse problems as well. The same issue showed up, for instance, in the field of groundwater transmissivity detection \cite{giudici-parra}.

This kind of questions also arises in the branch of the homogenization theory which studies effective properties of composite materials. We give a brief outline here.

Indeed, the positivity of Jacobians of injective $\sigma$-harmonic mappings has attracted attention in several applications. In two dimensions, the first application of this positivity has been given in \cite{nesi}. The long standing problem of improving the so-called Hashin-Shtrikman bounds \cite{hashin-shtrikman} for the effective conductivity of composite materials was addressed in that paper. The method used was based on ideas of Murat and Tartar \cite{murat-tartar} and Tartar  \cite{tartar79}, a reference not easy to find. We refer to \cite{tartar} for a more complete treatment. 

\noindent
The bottom line is as follows. The question of interest, in the simplest not yet known at that time, case is the following. Three  numbers $0<\sigma_1<\sigma_2<\sigma_3$, representing the conductivity of three isotropic materials, called the phases,  three ``volume fractions'' $p_1,p_2$ and $p_3$, summing up to 1 and representing the area proportions of the phases and  a $2\times 2$ matrix $A$, parametrizing the affine boundary data, are given. Assume that $\sigma=\displaystyle{\sum_{i=1}^3 \chi_i(x) \sigma_i}$ where $\chi_i$ represents the characteristic function of the set where $\sigma$ is equal to $\sigma_i$ times the identity matrix and $\frac{1}{|\Omega|} \int_{\Omega} \chi_i(x) dx = p_i, i=1,2,3$. Then one aims to determine a bound from below for the following quantity
\begin{equation}\label{vp}
F(A)=\inf_{\chi_1,\chi_2,\chi_3} \inf_{U_0\in W^{1,2}_0(\Omega)}\frac{1}{|\Omega|} \int_{\Omega} {\rm Trace}[(DU_0(x)+A)^T \sigma(x) (DU_0(x)+A)] dx\,.
\end{equation}
The overall problem is non linear and actually is linked with the notion of quasi-convexity. More precisely, it computes the quasiconvexification at the matrix $A$, of a non-convex function of $DU$, with $U-A x\in W^{1,2}_0(\Omega)$. This function turns out to be the minimum of three quadratic functions, as shown by Kohn and Strang in  \cite{kohn-strang}. However, for our purposes, it is important to note that the infimum over $U_0$ in \eqref{vp} is attained exactly when $U(x)=U_0(x)+Ax$ is the  $\sigma$-harmonic mapping with affine boundary data given by $U=A x$ on $\partial \Omega$. Optimal lower bound for \eqref{vp} were found by Kohn and Strang exactly exploiting the connection with the optimal bound for effective conductivity found by Murat and Tartar and, later, by Cherkaev and Lurie \cite{lurie-cherkaev}. The optimality is restricted to the case when only two isotropic phases are present that is,  only two materials are ``mixed''. For three or more phases, the methods based on compensated compactness gave suboptimal bounds. In this specific context, the compensated compactness method uses simply the constraint that the Jacobian determinant of the matrix $DU$  is a null-lagrangian.
 
\noindent
The classic strategy gives the so-called ``Wiener bound'', that is the harmonic mean bound. It is obtained considering the test fields $B$ in the class
\begin{equation*}\mathcal{B}_0:=\left\{B\in L^2(\Omega): \frac{1}{|\Omega|} \int_{\Omega} B(x)dx=A\right\}.
\end{equation*}
One obtains
\begin{equation*}
F(A)\geq  F_0(A):=\inf_{B\in \mathcal{B}_0}
 \frac{1}{|\Omega|} \int_{\Omega} {\rm Trace}[B(x)^T \sigma(x) B(x)] dx=
\end{equation*}
\begin{equation}\label{Wiener}
{\rm Trace}\left[A^T  \left(\frac{1}{|\Omega|} \int_{\Omega} \sigma^{-1}(x) dx\right)^{-1} A\right].
\end{equation}
\noindent
Tartar's ideas, based on compensated compactness, in this simplified context lead to an improved bound (called the ``translation bound'' by G. W. Milton) obtained by considering the new test field in the class
\begin{equation*}
\mathcal{B}_1:=\left\{B\in L^2(\Omega): \frac{1}{|\Omega|} \int_{\Omega} B(x)dx=A,  
 \frac{1}{|\Omega|} \int_{\Omega} \det B(x)dx=\det A\right\}\,. 
 \end{equation*}
 One obtains 

\begin{equation*}
F(A)\geq   F_1(A):=\inf_{B\in \mathcal{B}_1}
 \frac{1}{|\Omega|} \int_{\Omega} {\rm Trace}[B(x)^T \sigma(x) B(x)] dx\,.
\end{equation*}
The computation is more involved than \eqref{Wiener}.  This technique, however,  gives an optimal answer in two dimensions, when specialized to the case of two-phase isotropic materials. When one deals with more than two phases this approach is no longer optimal. The results in \cite{BMN} have the following corollary.
Set 
\begin{equation*}
\mathcal{B}_2:=\left\{B\in L^2(\Omega): \frac{1}{|\Omega|} \int_{\Omega} B(x)dx=A,  
 \frac{1}{|\Omega|} \int_{\Omega} \det B(x)dx=\det A \det B\geq 0\,, \hbox{a.e. in $\Omega$} \right\}.
 \end{equation*}
  One has
 \begin{equation}\label{vpe}
F(A)\geq   F_2(A):=\inf_{B\in \mathcal{B}_2}
 \frac{1}{|\Omega|} \int_{\Omega} {\rm Trace}[B(x)^T \sigma(x) B(x)] dx.
\end{equation}
In fact, in  \cite{nesi}, it is proved that $F_2(A)>F_1(A)$, as soon as one deals with more that two isotropic phases, for suitable choices of the given parameters $p_i$ and matrices $A$. Later, new optimal microgeometries  were found for multiphase materials in \cite{achn} and, using again the positivity of the Jacobian determinant, it was possible to prove their optimality according to a stricter criterion, see \cite{acn}.
The key is exactly the universal bound given on the Jacobian determinant, which, in this context reads as the inequality $\det A \det B\geq 0$ in \eqref{vpe}.
In this context, it is highly desirable not to have any constraint on the regularity of the interfaces between phases. When $\sigma$ is non-symmetric, applications to composites have been given, for instance, in the context of the classic Hall effect by Briane and Milton \cite{briane-milton1, briane-milton2}. Other applications have considered the problem of determining which electric fields are realizable by Briane, Milton, and Treibergs\cite{bmt}.
On the other hand one would like to have similar improvements in higher dimensions. Briane and Nesi \cite{briane-nesi} studied the case of laminates of high rank showing that, for these special microgeometries the positivity of the determinant of the ``Jacobians'' of the corrector matrix holds in any dimension. To explain the result in detail would require too long a digression. However, roughly speaking, one could expect that in higher dimensions, even for discontinuous $\sigma$ one could hope for the positivity of the Jacobian determinant if one makes assumptions on the ``microgeometry''. On the other hand, even in the very restricted setting of periodic boundary conditions, particularly adapted to composites, and even under the assumption of dealing with only two isotropic phases, there is no hope to control the sign of the Jacobian determinant of $\sigma$-harmonic mappings without further assumptions on the nature of the interfaces. One explicit example was provided by Briane, Milton and Nesi \cite{briane-milton-nesi}.

\noindent
We now go back to the precise subject of the present paper.
We pose the following problem. 

\

\noindent
{\bf Problem 1.}

\noindent
Can we find Dirichlet data
\begin{equation}\label{D-data}
\Phi=(\phi_1,\cdots,\phi_n):\partial \Omega \to \mathbb R^n
\end{equation}
such that the corresponding solution mapping $U=(u_1,\cdots,u_n)$ is such that $\det DU$ is bounded away from zero \emph{independently} of the conductivity $\sigma$?

\

\noindent
Note that, in this context, it is essential that the choice of the boundary data is independent of $\sigma$, because $\sigma$ is the real unknown of the original inverse problem. As is easily understandable, some a-priori assumptions on $\sigma$, such as ellipticity, and some kind of regularity shall be needed.

\noindent


Problem 1 has a different phenomenology depending on the space dimension. When $n=2$ the issue is more or less completely understood, whereas when $n=3$ or higher, various kinds of pathologies show up. A review of such pathologies and a discussion of the open issues when $n\geq 3$ shall be the object of Section \ref{3d}.

\noindent
The principal aim of this note is to provide, when $n=2$, a quantitative lower bound on the Jacobian determinant under essentially minimal regularity assumptions. This is the content of our main Theorem \ref{theorem2} which is the new contribution of this paper to this subject.

We start reviewing the main known results in dimension $n=2$. It was proved in Bauman et al \cite{BMN} that, if $\sigma$ is H\"older continuous, $\Omega$ has $C^{1,\alpha}$ boundary and $\Phi$ is a $C^{1,\alpha}$ diffeomorphism onto the boundary of a convex domain, then $\det DU>0$ everywhere. Note that in \cite{BMN}, only symmetric matrices $\sigma$ were explicitly considered, however, in view of classical results on two dimensional elliptic first order systems with H\"older coefficients see, for instance, \cite{bersni} Appendix and also \cite{bojd'on} Proposition 5.1, the result extends as well to the non--symmetric case.
On the other hand, the present authors \cite{AN-ARMA}, proved that when $\sigma$ is merely $L^{\infty}$ and $\Phi$ is a homeomorphism onto the boundary of a convex domain, then $\det DU>0$ almost everywhere. In fact it was proved that, for every locally invertible, sense preserving, $\sigma$-harmonic mapping $U$ one has 
\begin{equation}\label{BMO}
\log \det DU \in BMO
\end{equation}
and, subsequently  \cite{AN-AASF}, this result was improved to
\begin{equation}\label{Pesp}
\det DU \in A_{\infty}
\end{equation}
that is the class of Muckenhoupt weights \cite{coifeff}.

\noindent
We recall that for purely harmonic mappings, Lewy's Theorem \cite{lewy}, states that for two-dimensional harmonic homeomorphisms, the Jacobian determinant cannot vanish at interior points. Hence, when $n=2$, harmonic homeomorphisms are, indeed, diffeomorphisms.  However the Jacobian determinant may vanish at boundary points.

It is also worth mentioning that the convexity assumption on the target of the boundary mapping $\Phi$ is sharp, Choquet \cite{Choquet}, Alessandrini and Nesi \cite{AN-Pisa}, if one wishes to have a condition expressed merely on the ``shape'' of the target and not on its parametrization.

Conversely, note that when no regularity is assumed on $\sigma$, the essential infimum of $\det DU$ on compact subsets of $\Omega$ might indeed be zero. In Section \ref{meyers}, an example, based on a well-known one by Meyers, is illustrated. 

\noindent
In the next Section \ref{core} we shall prove a quantitative version of the  result in \cite{BMN}.
The starting point relies on prescribing some quantitative assumption on the boundary data $\Phi$, when viewed as a parametrization of the boundary of the convex target, see Definitions \ref{Qfi}, \ref{gammaqc}, \ref{bdrymap}. The subsequent step consists on a quantitative lower bound of the modulus of the gradient of a scalar solution to equation \eqref{basicpde}, Theorem \ref{theorem1}. This estimate may be interesting on its own. Finally we state and prove our main result, Theorem \ref{theorem2}.
\section{The quantitative bounds.}\label{core}
\hspace{0.7 cm}

\noindent
Let $\phi:\mathbb R\to \mathbb R$ be a $T$-periodic $C^1$ function. Let $\omega:[0,\infty)\to[0,\infty)$ be a continuous strictly increasing function such that $\omega(0)=0$.
\begin{definition}\label{Qfi}
Given $m, M\in \mathbb R$, $m<M$, we say that $\phi$ is quantitatively unimodal 
if there exists numbers $t_1\leq t_2< t_3\leq t_4<t_1+T$ such that
\begin{equation}
\begin{array}{rllllll}
\phi(t)&=&m&t\in[t_1,t_2]\,,\\
\phi(t)&=&M&t\in[t_3,t_4]\,,\\
\phi^{\prime}(t)&\geq&\min\{\omega(t-t_2),\omega(t_3-t)\},&t\in[t_2,t_3]\,,\\
-\phi^{\prime}(t)&\geq&\min\{\omega(t-t_4),\omega(t_1+T-t)\},&t\in[t_4,t_1+T]\,.
\end{array}
\end{equation}
In the sequel we will refer to the quadruple $\{T,m,M,\omega\}$ as to the ``character of unimodality'' of $\phi$.
\end{definition}
\noindent
The concept of unimodality, but not this terminology, first appears in Kneser \cite{Kneser}, when he proved Rad\`o's conjecture \cite{rado} concerning the case of  ``purely'' harmonic mappings. The terminology ``unimodality'' was introduced in this context by Leonetti and Nesi \cite{Leonetti-Nesi}, following the work of Alessandrini and Magnanini \cite{AMsiam}. A different terminology (almost two-to-one functions) has also been used for the same concept, Nachman, Tamasan and Timonov\cite{nachman-tamasan}.

\noindent
Let $\Gamma\subset \mathbb R^2$ be a simple closed curve parametrized by a  $T$-periodic $C^1$ mapping
\begin{equation}
\Phi: \mathbb R\to \mathbb R^2
\end{equation}
in such a way that $\Phi_{|_{[0,T)}}$ is one-to-one.
\begin{definition}\label{gammaqc}
We say that $\Gamma$ is quantitatively convex 
if for every $\xi\in \mathbb R^2$, $|\xi|=1$ the function
\begin{equation*}
\phi_{\xi}=\Phi\cdot \xi
\end{equation*}
is quantitatively unimodal and its character of unimodality is given by $\{T,m_{\xi},M_{\xi},\omega \}$ with $m_{\xi},M_{\xi}$ such that $M_{\xi}-m_{\xi}\geq D$, for a given $D>0$.

\noindent
In the sequel we will refer to the triple $\{T,D,\omega\}$ as to the ``character of convexity'' of $\Gamma$.
\end{definition}
\begin{remark}
If $\Gamma$ is quantitatively convex then it is convex, that is, it is the boundary of a convex set $G$. In fact each tangent line to $\Gamma$ turns out to be a support line for $G$. The following Lemma  provides a sufficient condition for quantitative convexity. 
Roughly speaking, it says that if $\Gamma$ is an appropriately  parametrized $C^2$ simple closed curve with strictly positive curvature, then it is quantitatively convex in the sense of Definition \ref{gammaqc}, and the character of  convexity can be computed in terms of the parametrization. Here, for the sake of simplicity, we have chosen the arc--length parametrization, because the main purpose of this Lemma is to provide a variety of examples, but we emphasize that in general, the character of convexity does depend on the parametrization of the curve and not only on its image.
\end{remark}

\noindent
We convene to denote by  $J$ the matrix representing the counterclockwise rotation of $90$ degrees
$$
J=
\left(
\begin{array}{cc}
0&-1\\
1&0
\end{array}
\right)\,.
$$
\begin{lemma}\label{lemmaQc}
Let $\Gamma$ be such that $\Phi\in C^2$ and:
\begin{equation}\label{lemma2}
\begin{array}{lll}
i)&|\Phi^{\prime}|=1\,,\\
ii)&0<\kappa\leq \Phi^{\prime\prime}\cdot J^T\Phi^{\prime}\leq K\,,
\end{array}
\end{equation}
 then $\Gamma$ is quantitatively convex with character $\{|\Gamma|,\frac 1 K,\frac{2\kappa}{\pi}\,t \}.$
\end{lemma}
\proof 

\noindent Condition $i)$ of Lemma \ref{lemmaQc}, implies that $\Phi^{\prime}(t)= e^{i s(t)},\, 0\leq t\leq T$ and we may assume that $s(0)=0.$ 
Without loss of generality we assume that $\Phi$ is orientation preserving. Then, by condition $ii)$ of Lemma \ref{lemmaQc}, one has $0<\kappa\leq s^{\prime}(t)\leq K.$ Picking, w.l.o.g., $\xi=e_2$,
\begin{equation*}
\phi_{\xi}(t)-\phi_{\xi}(0)=\int_0^t \sin(s(\tau)) d\,\tau\,.
\end{equation*}
The function $s(t)$ ranges over the whole interval $[0,2\pi]$, picking $t_{\pi}$ such that $s(t_{\pi})=\pi$, we have
\begin{equation*}
\begin{array}{lll}
M_{\xi}-m_{\xi}=\int_0^{t_{\pi}} \frac{\sin(s(\tau))}{
s^{\prime}(\tau)
}
\,ds(\tau)\,\geq \frac 1 K \int_0^\pi \sin(s) ds =\frac 1 K\,,\\
\\
\phi_{\xi}(t)=\sin(s(t))\,,\\
\\
s(t)= \int_0^t s^{\prime}(\tau)d\tau\geq\kappa\,t\,,\\
\\
\phi_{\xi}(t)\geq \frac{2\kappa}{\pi}t \,,\hskip0,2cm0\leq s(t)\leq \frac \pi 2\,.
\end{array}
\end{equation*}
Thus we may pick $D=\frac 1 K$ and $\omega(t)= \frac{2\kappa}{\pi}\,t\,\,,t\geq0\,.
$
\endproof

\noindent
We shall consider $\Omega$ a bounded simply connected domain in $\mathbb R^2$ with $C^{1,\alpha}$ boundary. In order to make precise the quantitative character of such regularity we introduce the following definition.

\begin{definition}\label{Qdomain}
A domain $\Omega\subset \mathbb R^2$ is said to be of class $C^{1,\alpha}$ with constants $\rho_0, M_0$,  positive and H\"older exponent $\alpha\in(0,1]$, if for any $P\in \partial \Omega$, there exist a rigid change of coordinates such that $P=0$ and we have
\begin{equation}
\Omega\cap B_{\rho_0}(0)=\{x\in B_{\rho_0}(0): x_2>\psi(x_1)\},
\end{equation}
where $\psi:[-\rho_0,\rho_0]\to \mathbb R^2$ is a $C^{1,\alpha}$ function satisfying 
\begin{equation}
\psi(0)=\psi^{\prime}(0)=0
\end{equation}
and also
\begin{equation}
||\psi||_{L^{\infty}([-\rho_0,\rho_0])}+\rho_0 ||\psi^{\prime}||_{L^{\infty}([-\rho_0,\rho_0])}+\rho_0^{1+\alpha}
\sup_{\stackrel{x,x^{\prime}\in [-\rho_0,\rho_0]}{ x\neq x^{\prime}}}\frac{|\psi^{\prime}(x)-\psi^{\prime}(x^{\prime})|}{|x-x^{\prime}|^{\alpha}}\leq M_0\rho_0\,.
\end{equation}
\end{definition}

\noindent
\begin{definition}\label{bdrymap} Given a $C^{1,\alpha}(\partial \Omega;\mathbb R)$ function $\phi$, we shall say that it is quantitatively unimodal, if considering the arclength parametrization of $\partial \Omega$, $x=x(s)$, $0\leq s\leq T=|\partial \Omega|$, the periodic extension of the function  $[0,T]\ni s\to \phi(x(s))$ is quantitatively unimodal with character $\{T,m,M,\omega\}$. For such a function $\phi$, we introduce the following closed arcs, possibly collapsing to a single point:
\begin{equation}
\begin{array}{lll}
\Gamma_{\rm min}=\{x\in \partial \Omega: \phi=m\}\,,\\
\Gamma_{\rm max}=\{x\in \partial \Omega: \phi=M\}\,.\\
\end{array}
\end{equation}
Accordingly, a mapping $\Phi \in C^{1,\alpha}(\partial \Omega;\mathbb R^2)$ shall be said quantitatively convex with character $\{T,D,\omega\}$ if the periodic extension of $\Phi(x(s))$ fulfils the conditions of Definition \ref{gammaqc}.
\end{definition}

\noindent
Let us consider $\sigma=\{\sigma_{ij}\}_{i,j=1,2}$ a, not necessarily symmetric, matrix of coefficients
$\sigma_{ij}:\overline{\Omega}\to\mathbb R$ satisfying the ellipticity condition
\begin{equation}\label{ell}
\begin{array}{ccrllll}
\sigma(x) \xi\cdot \xi&\geq& K^{-1} |\xi|^2&,&\hbox{for every $\xi \in \mathbb R^2$}\,,\\
\sigma^{-1}(x) \xi \cdot \xi &\geq & K^{-1}
|\xi|^2&,&\hbox{for every $\xi \in \mathbb R^2$}
\end{array}
\end{equation}
for given positive constant $K$,
and also
\begin{equation}\label{Hell}
\begin{array}{ll}
|\sigma_{ij}(x)-\sigma_{ij}(x^{\prime})|\leq E |x-x^{\prime}|^{\alpha}\,,&\forall x,x^{\prime}\in \overline{\Omega}\, ,
\end{array}
\end{equation}
for given $\alpha$, $0<\alpha\leq 1$ and $E>0$\,.

\noindent
We shall consider the $W^{1,2}(\Omega)$ solution $u$ to the Dirichlet problem
\begin{equation}\label{dp}
\left\{
\begin{array}{lll}
{\rm div} (\sigma \nabla u)=0&{\rm in}&\Omega\,,\\
u=\phi&{\rm on}&\partial \Omega\,.
\end{array}
\right.
\end{equation}
We recall that, in view of the classical regularity theory, $u$ in fact belongs to $C^{1,\beta}(\overline{\Omega})$, for some $\beta\leq \alpha$ and its norm is dominated by the $C^{1,\alpha}$-norm of $\phi$, modulo a constant which only depends on $\rho_0, M_0, K$ and $E$, with $\rho_0, M_0$ as in Definition \ref{Qdomain}.

\begin{lemma}\label{lemma1}
Let $\phi:\partial \Omega\to \mathbb R$ be quantitatively unimodal with character $\{|\partial \Omega|,m,M,\omega\}$ and assume that
\begin{equation}\label{Hfi}
\begin{array}{ll}
\left|
\frac{d}{ds}\phi(x(s))-\frac{d}{ds}\phi(x(s^{\prime}))
\right|
\leq E |s-s^{\prime}|^{\alpha}\,,&\forall s, s^{\prime}\in [0,|\partial \Omega|]\,.
\end{array}
\end{equation}
Then there exist $\kappa, \delta$ only depending on the character of unimodality (see Definitions \ref{Qfi}, \ref{bdrymap}) and on $\alpha, E$, such that 
if
\begin{equation}
\begin{array}{llllll}
x\in \overline{\Omega}&{\rm and}& {\rm dist}(x, \Gamma_{\rm min}\cup \Gamma_{\rm max})\leq \delta\,,
\end{array}
\end{equation}
then
\begin{equation}
|\nabla u(x)|\geq \kappa.
\end{equation}
\end{lemma}
\proof

\noindent
 Up to a $C^{1,\alpha}$ diffeomorphism, with constants only depending on $\rho_0, M_0$ and $|\partial \Omega|$, we may assume that $\Omega=B_1(0)$.

\noindent
It is well known that in such new coordinates $u$ solves a new Dirichlet problem of type \eqref{dp} with a new matrix of coefficients and  new boundary data that, however, satisfy analogous assumptions with constants and parameters only depending on the same a-priori data. For the sake of not to overburn the notation we stick to the  one of \eqref{dp}.

\noindent
By the $C^{1,\beta}$ regularity of $u$, if ${\rm dist}(x,\Gamma_{\rm max})\le \eta$, then $u(x)\geq M- C\eta$ with $C>0$ only depending on the a-priori data.

\noindent
Let us pick $\eta$ such that
\begin{equation*}
M-C\eta\geq \frac{M-m}{2}\,.
\end{equation*}
Hence, by Harnack's inequality \cite{G-T},
\begin{equation*}
\begin{array}{lll}
u(x)-m\geq C \, \frac{M-m}{2}>0,&\hbox{for every}&x\in B_{1-\eta}(0)\,.
\end{array}
\end{equation*}
\noindent
Here $C$ only depends on the a-priori data. By the version of the Hopf Lemma due to Finn and Gilbarg \cite{F-G} Lemma 7, which applies to equations in divergence form, and H\"older continuous $\sigma$, we obtain
\begin{equation}
\begin{array}{ll}
|\nabla u(x)|\geq \kappa_0>0\,,&\forall x\in \Gamma_{\rm min}\,,
\end{array}
\end{equation}
with $\kappa_0$ only depending on the a-priori data.

\noindent
By $C^{1,\beta}$ regularity we have
\begin{equation*}
\begin{array}{llllll}
|\nabla u(x)|\geq \kappa_0-C\delta^{\beta}\,&\forall x\in \overline{\Omega}&\hbox{such that}&{\rm dist}(x,\Gamma_{\rm min})\leq \delta\,.
\end{array}
\end{equation*}
Picking $\delta$ such that $C\delta^{\beta}\leq \frac{\kappa_0}{2}$, we obtain
\begin{equation}
\begin{array}{lllll}
|\nabla u(x)|\geq \frac{\kappa_0}{2}>0\,,&\hbox{if}&{\rm dist}(x,\Gamma_{\rm min})\leq \delta\,.
\end{array}
\end{equation}
A symmetrical result applies in the neighborhood of $\Gamma_{\rm max}$.
\endproof
\begin{lemma}\label{lemma3}
Under the same assumptions as in Lemma \ref{lemma1}, there exists $r>0$ such that
\begin{equation}
\begin{array}{lllll}
|\nabla u(x)|\geq L>0\,,&\forall x\in \overline{\Omega}\,,&{\rm dist}(x,\partial \Omega)\leq  r\,.
\end{array}
\end{equation}
Here $L$ and $r$ are positive and only depend on the a-priori data.
\end{lemma}

\proof

\noindent If we pick $x\in \partial \Omega$, and write $x=x(s)$ such that ${\rm dist}(x,\Gamma_{\rm min}\cup \Gamma_{\rm max})\geq \delta$, we have
\begin{equation*}
|\nabla u(x(s)\cdot x^{\prime}(s)|=\left|\frac{d}{ds}\phi(x(s)\right|\geq \omega(\delta)\,.
\end{equation*}
By $C^{1,\beta}$ regularity
\begin{equation*}
\begin{array}{lllll}
|\nabla u(x)|\geq \min\{\kappa,\omega(\delta)\}-C r^{\beta}\,,&\forall x\in \overline{\Omega}&\hbox{such that}&{\rm dist}(x,\partial \Omega)<r\,.
\end{array}
\end{equation*}
Picking $r$ such that 
\begin{equation*}
C r^{\beta}=\frac 1 2 \min\{\kappa,\omega(\delta)\}\,,
\end{equation*}
the thesis follows.
\endproof
\begin{theorem}\label{theorem1}
Let $\Omega$ be a simply connected domain, $C^{1,\alpha}$-regular with constants $\{\rho_0,M_0\}$ (see Definition \ref{Qdomain}).
Let $\phi:\partial \Omega\to \mathbb R$ be quantitatively unimodal with given character $\{|\partial \Omega|,m,M,\omega\}$(see Definitions \ref{Qfi}, \ref{bdrymap}) and let it satisfy the H\"older condition \eqref{Hfi}.
Let $\sigma=\{\sigma_{ij}(x)\}_{i,j=1,2}$ satisfy the ellipticity condition \eqref{ell} and the H\"older bound \eqref{Hell}.
Let $u\in W^{1,2}(\Omega)$ be the solution of the Dirichlet problem \eqref{dp}.

\noindent
Then there exists $C>0$, only depending on the a-priori data as above, such that
\begin{equation}\label{LB}
\begin{array}{lll}
|\nabla u (x) |\geq C>0\,, &\hbox{for every}& x\in \overline{\Omega}\,.
\end{array}
\end{equation}
\end{theorem}
\begin{remark}
Under stronger regularity assumptions, in particular assuming that $\sigma$ is Lipschitz continuous, a similar result was proven already in \cite{A-Pisa}, Theorem 3.2.
\end{remark}

\proof

\noindent As is well-known, there exists $\tilde{u}\in W^{1,2}(\Omega)$, called the stream function associated to $u$, which satisfies
\begin{equation}\label{sf}
\begin{array}{llll}
\nabla \tilde{u}=J\sigma \nabla u&\hbox{everywhere in}&\Omega\,,
&J=
\left(
\begin{array}{cc}
0&-1\\
1&0
\end{array}
\right)\,.
\end{array}
\end{equation}
Using complex notation $z=x_1+ i x_2$, $f=u+ i \tilde{u}$, the system \eqref{sf} can be rewritten as
\begin{equation}\label{1storder}
\begin{array}{ll}
f_{\bar{z}}=\mu f_z +\nu \bar{f_z}\ & \hbox{in $\Omega$}\ ,
\end{array}
\end{equation}
where, the so called complex dilatations $\mu , \nu$ are given by
\begin{equation}\label{SNU}
\begin{array}{llll}
\mu=\frac{\sigma_{22}-\sigma_{11}-i(\sigma_{12}+\sigma_{21})}{1+{\rm
Tr\,}\sigma +\det \sigma}& \ ,&\nu =\frac{1-\det \sigma
+i(\sigma_{12}-\sigma_{21})}{1+{\rm Tr\,}\sigma +\det \sigma}\ ,
\end{array}
\end{equation}
and satisfy the following
ellipticity condition
\begin{equation}\label{ellQC}
|\mu|+|\nu|\leq \frac{K-1}{K+1} \,,
\end{equation}
 and, being $\sigma$ H\"older continuous, also $\mu$ and $\nu$ satisfy an analogous H\"older bound.

\noindent
In \cite{BMN}, it is proven that $f$ is a $C^{1,\beta}$ diffeomorphism of $\overline{\Omega}$ onto $\overline{f(\Omega)}$. The lower bound obtained in Lemma \ref{lemma3}, implies that, setting
\begin{equation*}
\Omega_r=\{x\in \Omega: {\rm dist}(x, \partial \Omega)>r\}\,,
\end{equation*}
$f:\overline{\Omega}\backslash\Omega_r\to \mathbb C$,
is a bilipschitz homeomorphism with constants only depending on the a-priori data. We have identified $\mathbb C$ with $\mathbb R^2$ in the canonical way.

\noindent
Hence $f(\Omega)$ is also a $C^{1,\beta}$ domain with constants controlled by the a-priori data. Note that also $|\partial (f(\Omega))|$ is controlled.

\noindent
Let us denote $g=f^{-1}(w)$, $w\in \mathbb C$. A straightforward calculation gives
\begin{equation}\label{inverseB}
g_{\overline{w}}=-\nu(g) g_{w}- \mu(g)\overline{g_{w}}\,.
\end{equation}
In other words $g$ satisfies a Beltrami equation whose coefficients satisfy uniform ellipticity and H\"older continuity, with constants only depending on the a-priori data.

\noindent
By standard interior regularity estimates, $g_{w}$ is bounded in $f(\Omega_r)$. Using \eqref{inverseB}, we have
\begin{equation}\label{upperbound}
\begin{array}{lll}
|g_w|^2-|g_{\overline{w}}|^2\leq C^2& \hbox{in}&f(\Omega_r)\,,
\end{array}
\end{equation}
which can be rewritten as 
\begin{equation*}
\begin{array}{lll}
|f_w|^2-|f_{\overline{w}}|^2\geq C^{-2}& \hbox{in}&\Omega_r\,,
\end{array}
\end{equation*}
which in turn implies
\begin{equation}
\begin{array}{lll}
|\nabla u|\geq C^{-1}& \hbox{in}&\Omega_r\,.
\end{array}
\end{equation}
Hence, in combination with Lemma \ref{lemma3}, the thesis follows.
\endproof

\begin{theorem}\label{theorem2}
Let $\Omega$ and $\sigma$ be as in Theorem \ref{theorem1}. Let $\Phi=(\phi_1,\phi_2):\partial \Omega \to \mathbb R^2$ be quantitatively convex, see Definitions \ref{gammaqc}, \ref{bdrymap}, with character $\{|\partial \Omega|,D,\omega\}$. Let $U=(u_1,u_2)\in W^{1,2}(\Omega;\mathbb R^2)$ solve
\begin{equation}\label{dpv}
\left\{
\begin{array}{lll}
{\rm div} (\sigma \nabla u_i)=0&{\rm in}&\Omega\,,\\
u_i=\phi_i&{\rm on}&\partial \Omega\,.
\end{array}
\right.
\end{equation}
There exists $C>0$ only depending on the a-priori data such that
\begin{equation*}
U:\overline{\Omega}\to \overline{U(\Omega)}\subset \mathbb R^2
\end{equation*}
is a $C^{1,\beta}$ diffeomorphism and
\begin{equation}\label{det>c}
\begin{array}{lll}
\det DU \geq C^2>0&\hbox{in}&\overline{\Omega}\,.
\end{array}
\end{equation}
\end{theorem}
\begin{remark} A similar result, under slightly more restrictive hypotheses, has been recently proved by G. S. Alberti \cite{alberti}. In fact the approach in \cite{alberti} is based on estimates in \cite{A-Pisa} and \cite{A-ARMA} which require the Lipschitz regularity of $\sigma$.
Conversely, under somewhat different regularity assumptions, quantitative upper bounds on the so-called dilatation of a $\sigma$-harmonic mapping $U$, that is the quotient
$$
\frac{{\rm Trace} (DU^T  DU)}{2 \det DU}\,,
$$
have been recently studied in \cite{ANroma}, Theorems 3.1, 3.4 .
\end{remark}
\proof

\noindent
 In \cite{BMN}, it is shown that $U$ is a orientation preserving diffeomorphism. The lower bound \eqref{det>c} remains to be proven. For any $\xi\in \mathbb R^2$, $|\xi|^2=1$, we may apply Theorem \ref{theorem1} to 
 \begin{equation*}
 u_{\xi}=U\cdot \xi
 \end{equation*}
 and obtain
\begin{equation}
\begin{array}{lll}
|(DU) \xi|=|\nabla u_{\xi}|\geq C>0&\hbox{in}&\overline{\Omega}\,.
\end{array}
\end{equation}
Or equivalently
\begin{equation}
\begin{array}{lll}
|DU^T DU \xi\cdot \xi|=|\nabla u_{\xi}|^2\geq C^2>0&\hbox{in}&\overline{\Omega}\,,
\end{array}
\end{equation}
that is the eigenvalues $\lambda_1(x)$ and $\lambda_2(x)$ of the symmetric matrix $DU^T(x) DU(x)$  are uniformly bounded from below:
\begin{equation*}
\begin{array}{lll}
\lambda_i(x)\geq C^2>0\,,&i=1,2\,,&\forall \,x\in\overline{\Omega}\,.
\end{array}
\end{equation*}
Therefore
\begin{equation}
\begin{array}{lll}
(\det DU)^2=\lambda_1(x)\lambda_2(x)\geq C^4>0\,,&\hbox{everywhere in} &\overline{\Omega}\,.
\end{array}
\end{equation}
Since $U$ is sense preserving, one has $\det DU\geq C^2>0$ everywhere in $\Omega$.
\endproof
\begin{remark}
Theorem \ref{theorem2}, has some feature in common with the results in \cite{AN-Pisa}. In the latter paper the authors consider harmonic mappings which are extensions of given Dirichlet data aiming for univalent solutions. A characterization is given, for the case when $U$ is a diffeomorphism up to the boundary of $\Omega$, 
 in terms of the value of the Jacobian determinant on the boundary, so implicitly imposing constraints on the parametrization of the boundary of the image.
 One may wonder whether an assumption just on the shape of the target may suffice. This is not the case even in the purely harmonic case. Indeed one may exhibit a sequence $U_n$ of sense preserving, injective, harmonic mappings of the unit disk onto itself, fixing $U_n(0)=0$, such that $\det DU_n(0)\to 0$ as $n\to +\infty$. The convergence holds uniformly on compact subsets of the unit disk. The limit harmonic mapping in not univalent. See \cite{duren}, Section 4.1.
 \end{remark}
\section{Discontinuous coefficients. An example.} \label{meyers}
We elaborate on a well-known example by Meyers \cite{Meyers}. See also Leonetti and Nesi \cite{Leonetti-Nesi} for an application in a related context.
For a fixed $\alpha>0$ we consider the symmetric matrix of coefficients
\begin{equation}
\sigma(x)=
\left(
\begin{array}{cc}
\frac{\alpha^{-1}x_1^2 +\alpha x_2^2}{x_1^2+x_2^2}&\frac{(\alpha^{-1}-\alpha)x_1 x_2}{x_1^2+x_2^2}\\
\\
\frac{(\alpha^{-1}-\alpha)x_1 x_2}{x_1^2+x_2^2}&\frac{\alpha x_1^2 +\alpha^{-1} x_2^2}{x_1^2+x_2^2}
\end{array}
\right)\,.
\end{equation}
Is is a straightforward matter to check that its entries belong to $L^{\infty}$ and that $\sigma$ has eigenvalues $\alpha$ and $\alpha^{-1}$. Therefore $\sigma$ satisfies the uniform ellipticity condition \eqref{ell} with ellipticity constant
\begin{equation*}
K=\max\{\alpha,\alpha^{-1}\}\,,
\end{equation*}
and $\sigma$ is discontinuous at $(0,0)$ (and only at $(0,0)$, when $\alpha\neq 1$).
Let us denote
\begin{equation*}
\begin{array}{l}
u_1(x)=|x|^{\alpha -1}x_1\,,\\
u_2(x)=|x|^{\alpha -1}x_2\,.
\end{array}
\end{equation*}
A direct calculation shows that $u_i\in W^{1,2}(B_1(0))$, $i=1,2$ and that they solve the Dirichlet problem
\begin{equation*}
\left\{
\begin{array}{llc}
{\rm div} (\sigma \nabla u_i)=0&{\rm in}&B_1(0)\,,\\
u_i=x_i&{\rm on}&\partial B_1(0)\,.
\end{array}
\right.
\end{equation*}
Note also that $f=u_1+i u_2$ is a quasiconformal mapping of $B_1(0)$ onto itself and it solves the Beltrami equation 
\begin{equation*}
f_{\overline{z}}=\frac{\alpha-1}{\alpha+1} \,\frac{z}{\overline{z}} \,f_z\,.
\end{equation*}
Setting $U=(u_1,u_2)$, we compute 
\begin{equation*}
\det DU =|f_z|^2-|f_{\overline{z}}|^2=\alpha |z|^{2(\alpha-1)}\,.
\end{equation*}
Therefore $\det DU$ vanishes at $(0,0)$ when $\alpha>1$, whereas, when $\alpha\in (0,1)$, it diverges as $z\to 0$.
%

%
%
%
%
%
%
%
%
\section{Mappings in higher dimensions. Examples and open problems.}\label{3d}
The interior lower bound on $\det DU$ obtained in Theorem \ref{theorem2}, has been achieved by methods which are intrinsically two-dimensional (the Beltrami equation). Only part of the result can be extended to higher dimensions.

\noindent
For instance, with minor adaptations of the method developed in the Section \ref{core}, one can argue as follows.

\noindent
Consider $\Omega\subset \mathbb R^n$, a bounded domain diffeomorphic to a ball of class $C^{1,\alpha}$ and with constants  $\rho_0,M_0$ defined with the obvious slight adaptations of Definition \ref{Qdomain}. 

\noindent
Let $\sigma=\{\sigma_{ij}\}_{i,j=1,2}$ be the matrix of coefficients and let it satisfy uniform ellipticity with constant $K$ as in \eqref{ell} and H\"older continuity like in \eqref{Hell}. 

\noindent
Let $G\subset \mathbb R^n$ be a convex body whose boundary $\Gamma$ is $C^2$ and having at each point principal curvatures bounded from below by $\kappa>0$.

\noindent
Let $\Phi=(\phi_1,\phi_2,\cdots,\phi_n):\partial \Omega \to \Gamma$ be an orientation preserving diffeomorphism such that
$\Phi,\Phi^{-1}$ are $C^{1,\alpha}$ with constant $E$.
\noindent
Let $U=(u_1,u_2,\cdots,u_n)\in W^{1,2}(\Omega;\mathbb R^n)$ be the weak solution to
\begin{equation*}
\left\{
\begin{array}{llc}
{\rm div} (\sigma \nabla u_i)=0&{\rm in}&\Omega\,,\\
u_i=\phi_i&{\rm on}&\Omega\,,\\
i=1,2\cdots,n\,.
\end{array}
\right.
\end{equation*}
Then, by the same arguments used in Section \ref{core}, we obtain.
\begin{theorem}
Under the above stated assumptions, there exists $\rho>0$ and $Q>0$ such that $U$ is a diffeomorphism of $\overline{\Omega}\backslash \Omega_{\rho}$ onto a neighborhood of $\Gamma$, within $\overline{G}$ and we have
\begin{equation*}
\begin{array}{lll}\det DU\geq Q&\hbox{in}&\overline{\Omega}\backslash \Omega_{\rho}\,.
\end{array}
\end{equation*}
\end{theorem}
\noindent
We omit the proof.

\noindent
When $n\geq 3$, there is no chance, under the kind of hypotheses stated above, to obtain a global lower bound on $\det DU$. Evidence comes from a sequence of counterexamples that have been produced in a wide time span. A first illuminating example goes back to Wood \cite{wood} and has the amazing feature of being totally explicit. Wood displayed the following harmonic polynomial mapping from $\mathbb R^3$ onto $\mathbb R^3$:
\begin{equation*}
U(x_1,x_2,x_3)=(u_1,u_2,u_3)=(x_1^3-3x_1 x_3^2+x_2 x_3,x_2-3x_1 x_3,x_3)
\end{equation*}
that is $U$ is a homeomorphism, but not a diffeomorphism because $\det DU=0$ on the plane $\{x_1=0\}$.

Later Melas \cite{melas} provided an example  of a three dimensional harmonic homeomorphism $U:B_1(0)\to B_1(0)$ such that $\det DU(0)=0$.  Subsequently, Laugesen \cite{laugesen}, showed that there exists homeomorphisms $\Phi:\partial B_1(0)\to \partial B_1(0)$ which are arbitrarily close to the identity in the sup-norm such that the mapping $U=(u_1,u_2,\cdots,u_n)$ solving

\begin{equation*}
\left\{
\begin{array}{llc}
\Delta u_i=0&{\rm in}&B_1(0)\,,\\
U=\Phi&{\rm on}&\partial B_1(0)\,,\\
i=1,2\cdots,n\,.
\end{array}
\right.
\end{equation*}
is such that $\det DU$ changes its sign somewhere inside $B_1(0)$. Such examples are especially striking because they show that in dimension $n\geq 3$ it seems difficult to find a universal rule to select Dirichlet data in such a way that the corresponding harmonic (or $\sigma$-harmonic mapping is invertible at a topological level (because it may reverse orientation!) and not only as a differentiable mapping.

\noindent
One may wonder whether changing the topology of the boundary data may help. In the periodic case, obviously the harmonic functions are linear and this might have left the hope that, for variable coefficients, the periodic case may be better that the generic Dirichlet problem.  However this is not the case. In \cite{briane-milton-nesi}, it was proved that, in dimension three, one can find a matrix $\sigma$ taking only two values, proportional to the identity matrix, and a periodic  arrangement with a smooth interface, but such that the corresponding solution $U$ of the cell problem also reverses the orientation. The Jacobian determinant changes its sign in the interior of the (unit) cube of periodicity.

\noindent
If, from the above examples, it seems that few chances are left of finding a universal criterion by which choosing Dirichlet data such that, for each $\sigma$ (although smooth) the corresponding $\sigma$-harmonic mapping $U$ has nondegenerate Jacobian, then a more reasonable goal would be to find a way to control, in term of the Dirichlet data, the set of points where the Jacobian may degenerate and possibly evaluate the vanishing rate at such points of degeneration.

This appears as a completely open problem, not at all easy as the following example by Jin and Kazdan  \cite{kazdan-jin} shows.
Let $a\in C^{\infty}(\mathbb R;\mathbb R)$ and set
\begin{equation}
\sigma(x)=
\left(
\begin{array}{ccc}
1&a(x_3)&0\\
a(x_3)&1&0\\
0&0&b(x_3)
\end{array}
\right)\,,
\end{equation}
with 
\begin{equation*}
\left\{
\begin{array}{lllll}
a(x_3)= 0&\hbox{for}& x_3\leq 0\,,\\
a(x_3)\in(0, a_0) &\hbox{for}&x_3> 0&\hbox{with}&a_0\in(0,1)\,,\\
b(x_3)=\frac{1}{1-a^2(x_3)}&\hbox{for}&x_3\in \mathbb R\,.
\end{array}
\right.
\end{equation*}

We set
\begin{equation}
U(x) = \left(x_1, x_2, -x_1 x_2 +\phi(x_3)\right) \,,
\end{equation}
where $\phi$ is chosen in such a way that
\begin{equation*}
\left\{
\begin{array}{llc}
(b\phi^{\prime})^{\prime}-2 a = 0\,,&x_3\in \mathbb R\,,\\
\phi(x_3)=0\,,&x_3<0\,.
\end{array}
\right.
\end{equation*}
It turns out that $\phi^{\prime}>0$ for $x_3>0$ and consequently
\begin{equation*}
\det DU =\left\{
\begin{array}{llc}
\phi^{\prime}>0\,,&\hbox{for}&x_3>0\,,\\
\phi^{\prime}=0\,,&\hbox{for}&x_3\leq 0\,.
\end{array}
\right.
\end{equation*}
This means that the Jacobian determinant of a $\sigma$-harmonic mapping does not fulfill the property of unique continuation (whereas this is the case for $|DU|^2={\rm Trace}(DU^T DU)$). Hence the evaluation of the zero set of $\det DU$ from boundary data might be troublesome.
\begin{remark}
The above example has some striking features. First of all note also that, letting $a_0 \searrow 0$, we can make $\sigma$ as close as we want to the identity matrix. Moreover $U$ converges, uniformly on each compact subset of $\mathbb R^3$, to the harmonic polynomial mapping $U_0(x)= \left(x_1, x_2, -x_1 x_2 \right)$.
\end{remark}

We conclude by noticing that a limiting case of the above construction yields an example with a discontinuous, two--phase,  $\sigma$ which is remarkable as  well.

As before we pose
\begin{equation}
\sigma(x)=
\left(
\begin{array}{ccc}
1&a(x_3)&0\\
a(x_3)&1&0\\
0&0&b(x_3)
\end{array}
\right)\,,
\end{equation}
where now 

\begin{equation*}
\left\{
\begin{array}{lllll}
a(x_3)= 0&\hbox{for}& x_3\leq 0\,,\\
a(x_3)=a_0 &\hbox{for}&x_3> 0&\hbox{with}&a_0\in(0,1)\,,\\
b(x_3)=\frac{1}{1-a^2(x_3)}&\hbox{for}&x_3\in \mathbb R\,.
\end{array}
\right.
\end{equation*}

That is $\sigma$ is piecewise constant, namely 
\begin{equation*}
\begin{array}{ccccccc}
\sigma(x)=
\left(
\begin{array}{ccccc}
1&0&0\\
0&1&0\\
0&0&1
\end{array}
\right)\,
&\hbox{when}& x_3 < 0 
&\hbox{and}&
\sigma(x)=
\left(
\begin{array}{ccc}
1&a_0&0\\
a_0&1&0\\
0&0&\frac{1}{1-a_0^2}
\end{array}
\right)\,&\hbox{when}& x_3 > 0\,.
\end{array}
\end{equation*}
Again we pose
\begin{equation}
U(x) = \left(x_1, x_2, -x_1 x_2 +\phi(x_3)\right) \,,
\end{equation}
where now $\phi$ is given by
\begin{equation*}
\left\{
\begin{array}{llc}
\phi(x_3)= a_0(1-a_0^2) \,x_3^2\,,&x_3>0\,,\\
\phi(x_3)=0\,,&x_3\le 0\,.
\end{array}
\right.
\end{equation*}
We obtain that $U$ is a $\sigma$-harmonic mapping with $C^{1,1}$ regularity and, analogously to the previous example, it satifies
\begin{equation*}
\det DU =\left\{
\begin{array}{llc}
2 a_0(1-a_0^2) \,x_3>0\,,&\hbox{for}&x_3>0\,,\\
0\,,&\hbox{for}&x_3\leq 0\,.
\end{array}
\right.
\end{equation*}

\acknowledgement{G.A. was supported by FRA2012 `Problemi Inversi', Universit\`a degli Studi di Trieste, V.N. was supported by PRIN Project 2010-2011 `Calcolo delle Variazioni'.}

\end{document}